\documentclass{amsart}

\usepackage{amsfonts}
\usepackage{amssymb}
\usepackage{amsmath}

\newcommand{\hilb}{\mathbb{H}}

\newtheorem{theorem}{Theorem}[section]

\newtheorem{fact}[theorem]{Fact}
\newtheorem{question}{Open question}

\begin{document}

\title{On the complexity of topological conjugacy of compact metrizable $G$-ambits}
\author{Burak Kaya}
\address{ Department of Mathematics, Middle East Technical University, 06800, \c{C}ankaya, Ankara, Turkey\\
}
\email{burakk@metu.edu.tr}

\begin{abstract}
In this note, we analyze the classification problem for compact metrizable $G$-ambits for a countable discrete group $G$ from the point of view of descriptive set theory. More precisely, we prove that the topological conjugacy relation on the standard Borel space of compact metrizable $G$-ambits is Borel for every countable discrete group $G$.
\end{abstract}

\maketitle

\section{Introduction}

Let $G$ be a topological group. A $G$-\textit{flow} is a triple $(X,G,\Gamma)$ where $X$ is a compact Hausdorff space and $\Gamma: G \times X \rightarrow X$ is a continuous action of $G$ on $X$. A $G$-\textit{ambit} is a $G$-flow with a distinguished point whose $G$-orbit is dense, i.e. a quadruple of the form $(X,G,\Gamma,\alpha)$ where $(X,G,\Gamma)$ is a $G$-flow, $\alpha \in X$ and $\overline{\{\Gamma(g,\alpha): g \in G\}}=X$. From now on, for notational simplicity, we shall denote the point $\Gamma(g,x)$ by $g \cdot x$.

Two $G$-ambits $(X,G,\Gamma,\alpha)$ and $(Y,G,\Gamma',\beta)$ are said to be \textit{topologically conjugate} if there exists a homeomorphism $\pi: X \rightarrow Y$ such that $\pi(\alpha)=\beta$ and $\pi(g \cdot x)=g \cdot \pi(x)$ for all $g \in G$ and $x \in X$. From the point of view of topological dynamics, topologically conjugate $G$-ambits are considered the same. Our aim in this note is to analyze the complexity of the topological conjugacy problem for compact metrizable $G$-ambits for a countable discrete group $G$.

Over the last couple decades, a mathematical framework to rigorously compare the complexity of classification problems has been developed and successfully applied to diverse areas of mathematics. The idea is to represent classes of mathematical structures as \textit{Polish spaces}, i.e. completely metrizable separable topological spaces, and apply the techniques of descriptive set theory to analyze classification problems considered as definable equivalence relations on these Polish spaces. For a detailed development of this framework, we refer the reader to \cite{Gao09}.

In recent years, there has been a growing literature on the descriptive set theoretic analysis of classification problems from topological dynamics. For example, see \cite{Clemens09}, \cite{Thomas13}, \cite{GaoJacksonSeward15}, \cite{SabokTsankov15}, \cite{GaoHill16} and \cite{Kaya17}. In this note, we extend this study to compact metrizable $G$-ambits and prove the following.

\begin{theorem}\label{maintheorem} Let $G$ be a countable discrete group. Then the topological conjugacy relation on the standard Borel space of compact metrizable $G$-ambits is Borel.
\end{theorem}

Intuitively speaking, our result shows that there exists a countable (possibly, transfinite) procedure that uses countable amount of information for determining whether or not two compact metrizable $G$-ambits are topologically conjugate. On the other hand, the proof uses a classical theorem of Suslin and does not provide us with such an explicit procedure.

This note is organized as follows. In Section 2, we provide the necessary background from descriptive set theory. In Section 3, we construct the standard Borel space of compact metrizable $G$-ambits for a countable discrete group $G$. In Section 4, we prove Theorem \ref{maintheorem}. In Section 5, we discuss some issues regarding our coding of compact metrizable $G$-ambits and conclude with some open questions.

\section{Background from descriptive set theory}

In this section, we shall recall by some basic definitions and facts from descriptive set theory. For a general background, we refer the reader to \cite{Kechris95}.

A measurable space $(X,\mathcal{B})$ is called a \textit{standard Borel space} if $\mathcal{B}$ is the Borel $\sigma$-algebra of some Polish topology on $X$. It is well-known that if $A \subseteq X$ is a Borel subset of a standard Borel space $(X,\mathcal{B})$, then $(A,\{A \cap B: B \in \mathcal{B}\})$ is also a standard Borel space \cite[Corollary 13.4]{Kechris95}.

A subset of a Polish space $(X,\tau)$ is said to be \textit{analytic} if it is the projection of a Borel set $B \subseteq X \times Y$, i.e. it is of the form
\[\{x \in X: \exists y \in Y\ B(x,y) \}\]
for some Polish space $(Y,\tau')$. A subset of a Polish space $(X,\tau)$ is said to be \textit{coanalytic} if its complement is analytic, i.e. it is of the form
\[\{x \in X: \forall y \in Y\ B(x,y) \}\]
for some Polish space $(Y,\tau')$ and Borel relation $B \subseteq X \times Y$. It is clear that all Borel sets are analytic. On the other hand, there exist analytic sets that are not Borel \cite[Theorem 14.2]{Kechris95}. One can even find striking examples of such sets among classification problems. For example, Foreman, Rudolph and Weiss \cite{ForWeiss11} showed that the conjugacy relation on the space of ergodic measure preserving transformations is an analytic non-Borel set. This result may be interpreted as a strong anti-classification result since it shows that no ``explicit procedure" exists for determining whether or not two such transformations are conjugate. The following classical theorem of Suslin gives the relationship between Borel sets and analytic sets.
\begin{fact}\label{suslin} Let $S$ be a subset of a Polish space $(X,\tau)$. Then $S$ is Borel if and only if $S$ is analytic and coanalytic.
\end{fact}

Finally, we turn our attention to a special class of Polish spaces that will be used throughout this note. Let $(X,\tau)$ be a Polish space with a compatible complete metric $d$. It is well-known that the set $\mathcal{K}(X)$ of non-empty compact subsets of $X$ endowed with the topology induced by the Hausdorff metric
\[\delta_d(C_1,C_2)=\max\{\max_{x \in C_1} d(x,C_2), \max_{y \in C_2} d(y,C_1)\]
is a Polish space \cite[\S 4.F]{Kechris95}. 

\section{The standard Borel space of compact metrizable $G$-ambits}

In this section, we shall construct the standard Borel space of compact metrizable $G$-ambits for a countable discrete group $G$. Our construction basically follows one of the several constructions of compact metrizable $\mathbb{Z}$-flows given by Beleznay and Foreman \cite{ForBel95}.

For the rest of this note, let $G$ be a fixed countable discrete group. Since the collection of compact metrizable $G$-ambits is a proper class, in order to construct the standard Borel space of such $G$-ambits, we need to ``code" each compact metrizable $G$-ambit as an element of a set.

Let $\hilb$ denote the Hilbert cube $[0,1]^{\mathbb{N}}$. Recall that any compact metric space is homeomorphic to a compact subset of $\hilb$ and that this homeomorphism can be canonically constructed via the metric and a dense subset \cite[Theorem 4.14]{Kechris95}. Thus we can assume without loss of generality that the underlying topological space of a compact metrizable $G$-flow is a compact subset of $\hilb$.

Let $(X,G,\Gamma,\alpha)$ be a $G$-ambit. For each $g \in G$, the (graph of the) homeomorphism $\Gamma(g,\cdot)$ given by $w \mapsto g \cdot w$ is a compact subset of $X^2$ and hence of $\hilb^2$. We shall code the $G$-ambit $(X,G,\Gamma,x)$ as an element of the Polish space $(\mathcal{K}(\hilb^2))^G \times \hilb$ by representing each homeomorphism $\Gamma(g, \cdot)$ as an element of $\mathcal{K}(\hilb^2)$ in the corresponding component.

Let $\mathcal{A}$ be the subset of $(\mathcal{K}(\hilb^2))^G \times \hilb$ consisting of elements $((\lambda_g)_{g \in G},\alpha)$ satisfying the following conditions.
\begin{itemize}
\item[a.] $\forall g \in G$ ``$\lambda_g$ is a homeomorphism between compact subsets of $\hilb^2$"
\item[b.] $\forall g \in G\ \ dom(\lambda_g)=ran(\lambda_g)$
\item[c.] $\forall x \in \hilb\ \ \forall g,h \in G\ \ x \in dom(\lambda_g) \leftrightarrow x \in dom(\lambda_h)$
\item[d.] $\forall x \in \hilb\ \ x \in dom(\lambda_{1_G}) \rightarrow \lambda_{1_G}(x)=x$
\item[e.] $\forall x \in \hilb\ \ \forall g,h \in G\ \ x \in dom(\lambda_g) \rightarrow \lambda_g(\lambda_h(x))=\lambda_{gh}(x)$
\item[f.] $\alpha \in dom(\lambda_{1_G})$
\item[g.] $\forall x \in \hilb\ \ \forall \epsilon \in \mathbb{Q}^+\ \exists g \in G\ \ x \in dom(\lambda_{1_G}) \rightarrow d_{\hilb}(\lambda_g(\alpha), x) < \epsilon$
\end{itemize}
where $d_{\hilb}$ is a compatible complete metric for $\hilb$. We claim that $\mathcal{A}$ is Borel.

That conditions [a] and [b] define Borel subsets of $(\mathcal{K}(\hilb^2))^G \times \hilb$ follows from \cite[Lemma 3.5]{ForBel95}. Observe that that, provided that [a] and [b] hold, the universal quantifications over $\hilb$ in other conditions can be replaced by universal quantifications over a fixed countable dense subset of $\hilb$, since $\lambda_g$ is continuous and $dom(\lambda_g)$ is closed for all $g \in G$. Moreover, the relations
\[\{(s,C) \in S \times \mathcal{K}(S): s \in C\} \text{ and } \{(s,C) \in S \times \mathcal{K}(S \times S): s \in \pi_1(C)\}\]
are Borel (indeed, closed) relations for any Polish space $S$, where $\pi_1$ is the projection onto the first component. It easily follows that $\mathcal{A}$ is Borel and hence is a standard Borel space on its own.

The standard Borel space $\mathcal{A}$ is the space of compact metrizable $G$-ambits. In the next section, we shall prove that the topological conjugacy relation on $\mathcal{A}$ is a Borel subset of $\mathcal{A} \times \mathcal{A}$.

\section{Proof of Theorem \ref{maintheorem}}
Let $\mathcal{H}$ be the Borel subset of $\mathcal{K}(\hilb^2)$ satisfying condition [a] in the previous section, i.e. $\mathcal{H}$ is the standard Borel space of homeomorphisms between compact subsets of $\hilb^2$. By definition, two $G$-ambits $((\lambda_g)_{g \in G},\alpha),((\theta_g)_{g \in G},\beta) \in \mathcal{A}$ are topologically conjugate if and only if there exists $f \in \mathcal{H}$
\begin{itemize}
\item $dom(f)=dom(\lambda_{1_G})\ \wedge\ ran(f)=dom(\theta_{1_G})$
\item $\forall x \in \hilb\ \forall g \in G\ \ (x \in dom(\lambda_{1_G}) \rightarrow f(\lambda_g(x))=\theta_g(f(x)))$
\item $(\alpha,\beta) \in f$
\end{itemize}
As was the case before, in order to check the first two conditions, it suffices to quantify over a fixed countable dense subset of $\hilb$. Consequently, the subset of $\mathcal{H} \times \mathcal{A} \times \mathcal{A}$ defined by these three conditions is Borel and hence the topological conjugacy relation on $\mathcal{A}$ is an analytic subset of $\mathcal{A} \times \mathcal{A}$. It remains to prove that this relation is also coanalytic.

Let $(X,G,\Gamma,\alpha)$ and $(Y,G,\Gamma',\beta)$ be compact metrizable $G$-ambits. Since the $G$-orbits of $\alpha$ and $\beta$ are dense in $X$ and $Y$ respectively, we may naively attempt to continuously extend the map $g \cdot \alpha \mapsto g \cdot \beta$ from the $G$-orbit of $\alpha$ to all of $X$. More precisely, consider the relation $f$ defined by
\[ f(x)= \displaystyle \lim_{i \rightarrow \infty} g_i \cdot \beta \text{ if } \displaystyle x = \lim_{i \rightarrow \infty} g_i \cdot \alpha \text{ for some sequence } (g_i)_{i \in \mathbb{N}} \in G^{\mathbb{N}}\]
It is easily seen that any topological conjugacy between these $G$-ambits, if exists, is unique and has to be of this form. Why does this naive approach not prove that all compact metrizable $G$-ambits are topologically conjugate?

One reason is that it may be possible for $(g_i \cdot \alpha)$ to converge while $(g_i \cdot \beta)$ diverges, in which case the domain of $f$ is not $X$. Another reason is that it may be possible for $(g_i \cdot \alpha)$ and $(h_i \cdot \alpha)$ to converge to the same point while $(g_i \cdot \beta)$ and $(h_i \cdot \beta)$ converge to different points, in which case $f$ is not well-defined.

A moment's thought shows that these are the only obstructions preventing $f$ from being a continuous function with domain $X$, i.e. $f$ is a continuous function if $(g_i \cdot \beta)$ converges whenever $(g_i \cdot \alpha)$ converges, and $(g_i \cdot \beta)$ and $(h_i \cdot \beta)$ converge to the same point whenever $(g_i \cdot \alpha)$ and $(h_i \cdot \alpha)$ converge to the same point. Arguing symmetrically, one can also express necessary and sufficient conditions for $f$ being one-to-one and onto in terms of the convergence of sequences of points. Provided that $f$ is a homeomorphism, it is straightforward to check that $f$ intertwines the action of $G$ and that $f(\alpha)=\beta$. These observations provide us with the following (slightly complicated) characterization of topological conjugacy.

Two compact metrizable $G$-ambits $(X,G,\Gamma,\alpha)$ and $(Y,G,\Gamma',\beta)$ are topologically conjugate if and only if for all sequences $(g_i),(h_i) \in G^{\mathbb{N}}$, we have that
\begin{itemize}
\item[i.] $(g_i \cdot \alpha)$ converges (equivalently, is Cauchy) if and only if $(g_i \cdot \beta)$ converges (equivalently, is Cauchy)
\item[ii.] If $(g_i \cdot \alpha)$ and $(h_i \cdot \alpha)$ converge (equivalently, are Cauchy), then\\
$\displaystyle \lim_{i \rightarrow \infty} g_i \cdot \alpha=\displaystyle \lim_{i \rightarrow \infty} h_i \cdot \alpha \longleftrightarrow \displaystyle \lim_{i \rightarrow \infty} g_i \cdot \beta=\displaystyle \lim_{i \rightarrow \infty} h_i \cdot \beta$
\end{itemize}
We remark that whether or not a sequence $(g_i \cdot \alpha)$ is Cauchy can be checked in a Borel way. More precisely, the relation
\[\{((g_i),((\lambda_g)_{g \in G},\alpha)):\ \forall \epsilon \in \mathbb{Q}^+\ \exists k\ \forall m,n \geq k\ d_{\hilb}(\lambda_{g_m}(\alpha),\lambda_{g_n}(\alpha))<\epsilon\}\]
is a Borel subset of $G^{\mathbb{N}} \times \mathcal{A}$ since the set defined by the quantifier-free part of the condition is Borel for each $\epsilon \in \mathbb{Q}^+$ and $m, n \in \mathbb{N}$. Similarly,  whether or not two Cauchy sequences $(g_i \cdot \alpha)$ and $(h_i \cdot \alpha)$ have the same limit can be expressed by a Borel condition. Consequently, conditions [i] and [ii] define a Borel subset of $G^{\mathbb{N}} \times G^{\mathbb{N}} \times \mathcal{A} \times \mathcal{A}$. It follows that the topological conjugacy relation on $\mathcal{A} \times \mathcal{A}$ is also coanalytic and hence is Borel by Fact \ref{suslin}.
\section{Concluding remarks}

We remark that the Borel structure of the space of compact metrizable $G$-ambits is crucial for the proof of Theorem \ref{maintheorem}. It may not be obvious at first sight whether or not the same proof would work if we used a different coding of compact metrizable $G$-ambits. For example, there are three different constructions offered for the space of compact metrizable $\mathbb{Z}$-flows in \cite{ForBel95}. However, all these different constructions are proven to be equivalent up to Borel functions.

In practice, whenever there are different natural codings of the same class of mathematical structures as standard Borel spaces, these codings turn out to be equivalent up to Borel functions. Consequently, the analysis of the complexity of classification problems is not effected by the choice of coding. This phenomenon is so frequently observed that it is proposed in \cite[\S 14.1]{Gao09} as an analogue of the Church-Turing thesis for the descriptive set theoretic analysis of classification problems.

The proof of Theorem \ref{maintheorem} does not explicitly give a Borel condition characterizing topological conjugacy of compact metrizable $G$-ambits. On the one hand, such
Borel conditions have been found for certain classes of $G$-ambits. For example, the author proved in \cite{Kaya17} that the collection of countable atomless $\mathbb{Z}$-syndetic Boolean subalgebras of $\mathcal{P}(\mathbb{Z})$ can be used as a complete invariant for topological conjugacy of pointed Cantor minimal systems, i.e. $\mathbb{Z}$-ambits whose underlying space is the Cantor space and all of whose orbits are dense. On the other hand, as far as the author knows, no Borel condition is known for arbitrary compact metrizable $G$-ambits. We thus pose the following question.
\begin{question} With respect to Borel reducibility, what is the complexity of the topological conjugacy relation on the space of compact metrizable $G$-ambits for various countable discrete groups $G$?
\end{question}

\bibliographystyle{amsalpha}

\end{document}